\documentclass[report]{IEEEtran}
\usepackage{textcomp}
\usepackage{graphicx}
\usepackage[listofformat=subparens]{subfig}
\usepackage{comment}
\usepackage{amsmath}
\usepackage{amssymb}
\usepackage{graphicx}
\usepackage{algorithm}
\usepackage{algorithmic}
\usepackage{eurosym}
\usepackage{bbold}
\usepackage{booktabs} % To thicken table lines
\usepackage{fancyvrb}
\usepackage{bm}
\usepackage{dsfont}
\renewcommand{\algorithmicrequire}{\textbf{Input:}}

\newcommand{\x}{\bm{x}}

\newcommand{\xtime}{x_{\mathrm{t}}}

\newcommand{\action}{\bm{u}}
\newcommand{\w}{\bm{w}}

\newcommand{\uphysk}{u_{k}^{\textrm{phys},i}}

\newcommand{\uphysl}{u_{l}^{\textrm{phys}}}

\newcommand{\xphysk}{\x_{\textrm{phys},k}}

\newcommand{\xobski}{\bm{x}_{k}^{\mathrm{obs,i}}}

\newcommand{\probw}{p_{\mathcal{W}}(\cdot|x)}

\newcommand{\Xtime}{X_{\mathrm{t}}}
\newcommand{\Xphys}{X_{\mathrm{phys}}}

\usepackage{color,soul}
\soulregister\cite7
\soulregister\ref7
\soulregister\pageref7

\begin{document}

\title{Trajectory Tracking with an Aggregation of Domestic Hot Water Heaters: Combining Model-Based and Model-Free Control in a Commercial Deployment}

\author{Mingxi~Liu\IEEEauthorrefmark{2},~\IEEEmembership{Member,~IEEE,}, Stef~Peeters\IEEEauthorrefmark{4}, Duncan~S.~Callaway\IEEEauthorrefmark{2},~\IEEEmembership{Member,~IEEE,} and   Bert~J.~Claessens\IEEEauthorrefmark{4}\vspace{-0.45cm} 

\thanks{\IEEEauthorrefmark{2}M. Liu and D. S. Callaway are with the Energy \& Resources Group at University of California, Berkeley, 310 Barrows Hall, Berkeley, CA, 94720, USA (\{mxliu,dcal\}@berkeley.edu)} 
\thanks{\IEEEauthorrefmark{4}B.~J. Claessens and S.~Peeters work in the research department of REstore N.V. (bert.claessens@restore.energy).}
}

\markboth{\normalfont{Submitted to} Transactions on Smart Grid}{}

\maketitle

\begin{abstract}
Scalable demand response of residential electric loads has been a timely research topic in recent years. The commercial coming of age or residential demand response requires a scalable control architecture that is both efficient and practical to use. This work presents such a strategy for domestic hot water heaters and present a commercial proof-of-concept deployment. The strategy combines state of the art in aggregate-and-dispatch with a novel dispatch strategy leveraging recent developments in reinforcement learning and is tested in a hardware-in-the-loop simulation environment. The results are promising and present how model-based and model-free control strategies can be merged to obtain a mature and commercially viable control strategy for residential demand response.  
\end{abstract}


\IEEEpeerreviewmaketitle

\vspace{-0.10cm}
\section{Introduction}
\label{Sec.introductionDDR}

\IEEEPARstart{D}{emand} response programs can focus on a variety of applications, including  energy arbitrage \cite{mathieu2013energy}, ancillary services \cite{VrettosBoilers,hao2015aggregate} and voltage control \cite{weckx2014multiagent}.  A common challenge in the deployment of these programs is that of developing scalable and practical controls.  This challenge arises when controlling a large cluster of residential flexibility assets, for example, domestic hot water heaters (DHWHs).  These loads are a high potential source of flexibility for demand response programs \cite{VrettosBoilers,SMARTBOILER}, driven by their their decentralized abundance \cite{KAZMI2018159}, considerable and efficient storage capacity \cite{ClaessensBRLBoiler} and its negligible inertia. Algorithms that tap into this resource must take into account a variety of factors, including: (1) inter-temporal energy constraints, (2) the intrinsic uncertain user behavior, (3) methods to update models and control strategies as new information is collected and (4) the computational challenges associated with managing thousands to millions of devices to perform system-scale services.

Three well-studied paradigms for demand response control algorithms are that of model-based and model-free control and transactive energy based control \cite{PowerMatcher}.  Model-based control strategies start from a model of the flexibility assets which constrain an optimization problem that is solved at fixed time intervals. The poor scalability properties inherent to this approach can be mitigated by e.g. decomposition techniques or by working with a \textit{bulk} model \cite{mathieu2013energy} that describes the flexibility of a large cluster of flexibility assets in an \textit{aggregate-and-dispatch} approach. The bulk model is used to determine a control action for a cluster of assets, typically following a Model Predictive Control (MPC) strategy. The resulting control actions is dispatched over the assets using simple heuristics such as ranking according to the State of Charge (SoC) \cite{VrettosBoilers}.

An alternative to a model-based solution is to determine a control policy directly from data observed through interacting with the system to be controlled using techniques from reinforcement learning \cite{busoniu2010reinforcement}. Model-free control has the advantage of being more scalable as no new model needs to be engineered for each new asset, furthermore it does not suffer from model bias. This comes at the cost of a longer convergence time as adding prior knowledge can be cumbersome. Model-free control can be used to learn a control policy both at cluster and device level \cite{RuelensBRLCluster,RuelensBRLDevice}.

A third coordination mechanism is that of transactive energy control \cite{transactive}.  In this pragmatic approach, coordination is performed through decentralized market-based interaction.  However, though it is scalable and intuitive, it lacks \textit{planning} functionality as the market-based interactions typically supports only a myopic flexibility representation.

This paper presents a coordination strategy that could be used for real power control applications such as energy arbitrage, frequency regulation or peak demand management. It combines strong aspects of aggregate-and-dispatch and reinforcement-learning and has connections to transactive energy control in the sense that it dispatches taking into account opportunity costs represented by the advantage function learned directly from observed data. The control mechanism we propose is scalable, self adaptive and improves with data becoming available.  We also demonstrate the functionality of the controller with a network of real water heaters.  In Section \ref{Sec:related_work} an overview of the related literature is provided and the contributions of this work are explained. Section \ref{ImplementationDDR} presents the implementation of model-based and model-free control strategies. Assessment via networked simulations is presented in Section \ref{sec_simulation}, followed by hardware-in-the-loop test results in Section \ref{Sec.experiment}. Finally, Section \ref{Sec.Conclusions} outlines the conclusions and discusses future research.

%In Section~\ref{Sec:related_work} an overview of the related literature is provided and the contributions of this work are explained. Section~\ref{motivationes} sketches the main motivation behind this paper. Following the approach presented in~\cite{RuelensBRLDevice}, in Section~\ref{Sec:problem_formulation} a Markov decision process formulation is provided. In Section~\ref{ImplementationDDR}, the implementation of the controller is detailed, while Section~\ref{Sec:MABRL} presents a quantitative and qualitative assessment of its performance. Finally, Section~\ref{Sec:CDDR} outlines the conclusions and discusses future research.

\vspace{-0.10cm}
\section{Related Work and Contribution}\label{Sec:related_work}
\indent
In this Section an overview is given of aggregate-and-dispatch strategies and reinforcement learning applied to residential assets. 
\subsection{Aggregate-and-dispatch}\label{Sec:aggdispatch}
Aggregate-and-dispatch is a control method extensively studied in recent literature, the concept being that an aggregate model is derived representing the dynamics of a cluster of assets to be controlled. This model is used in an optimization problem that determines the aggregate set-point for the entire cluster.
This set-point is projected onto device-level actions using scalable heuristics.

The rationale for following this procedure is mainly driven by practical considerations, i.e. a reduced modeling effort and an optimization problem of reduced dimensionality making the decision making scalable.  In \cite{mathieu2013stateestimation,Liu_TIE_2016,Liu_TPS_2015}, state bin models are used to describe the dynamics of a large cluster of thermostatically controlled loads (TCLs). The bins cluster the TCLs according to their position relative to its individual temperature constraints, and linear models describe  transition dynamics between  different bins.  Model predictive control (MPC) strategies can coordinate switching actions in different bins~\cite{koch2011modeling}.  In \cite{Callawayheterogeneous}, a low-order tank model \cite{Kai} is used, allowing for a tractable stochastic optimization \cite{Leterme2014Flexible}, and in  \cite{mathieu2013energy} a low-order model demonstrated better results for the application of energy arbitrage. An extension to the tank model presented by Iacovella \textit{et al.} \cite{tracers}, is to use a set of representative TCLs to describe the dynamics of the cluster. This allows to model dispatch heuristics in the central optimization problem and accounts for the heterogeneous nature of the cluster of TCLs. In \cite{GoreckiTracking} an approach is presented that generates a \textit{trackable} reduced order tank model from a set of individual models using techniques from robust optimization.

A different approach to mitigate the limited scalability of a centralized solution is by relying on techniques from distributed optimization \cite{Gatsis2012,ADMMsmartgrid,Liu_TIE_2016}. Here the centralized optimization problem is decomposed over sub-problems. To obtain a global optimum, interaction between the sub-problems is required. Although demonstrating good optimization performance, this paradigm requires a model for each asset to be controlled, necessitating cumbersome automated system identification. Furthermore, the total computational and communication cost can be considerable as several iterations are required to obtain convergence. As the scope of this paper is on residential demand response where the value of the flexibility has to outweigh the cost to deploy the control solution this work targets an aggregate-and-dispatch approach as presented in \cite{Callawayheterogeneous,tracers}.

\subsection{Model-free dispatch}\label{Sec:RLreview}
In the papers discussed in Section \ref{Sec:aggdispatch} the dispatch strategy typically uses simple heuristics to decompose the aggregated set-point onto device level decisions. For example in \cite{tracers,Callawayheterogeneous,VrettosBoilers,ClaessensBRLBoiler} this is based upon the SoC of the DHWH, as this requires limited local intelligence and modeling effort. This however does not make a distinction between assets based upon energy efficiency, opportunity cost related to e.g. a reduced local consumption of renewable energy or the impact on the availability for demand response events. Taking opportunity costs into account when dispatching demand response events in a model-free way is one of the main contributions of this work.   This is achieved by leveraging recent results in data-driven control and more specifically reinforcement learning \cite{busoniu2010reinforcement}, requiring no explicit system identification. In order to do this one needs to know the opportunity costs for deviating from a control policy optimized for a local objective such as energy minimization. This can be obtained from the state-action value function or $Q$-function that can be learned directly from interactions of the controller with the DHWH.      
For example in \cite{ClaessensBRLBoiler,ClaessensHVAC,RuelensBRLDevice,oNeill2010residential,Hussein}, RL is applied in a residential setting to learn a control policy minimizing the cost of using energy. In \cite{ClaessensBRLBoiler,ClaessensHVAC,RuelensBRLDevice,oNeill2010residential} this is done by learning a state-action value function or Q-function which can be used directly to derive a control policy as detailed in Section \ref{ImplementationDDR}. The Q$^{\pi}$-function represents the value of being in a state $\x$ and taking action $\action$ following a policy $\pi$. This can be used to calculate an opportunity cost related to a demand response event that deviates from the local policy $h$ minimizing the cost for energy use and maximizing the availability for demand response events. An interpretation of this is that for each DHWH one uses RL to calculate a bid-function representing the cost for an action, this can be used in a market-based dispatch strategy such as transactive energy \cite{transactive,Vandael2012TSA} to dispatch the aggregated set-point.    
To summarize, the main contributions of this paper are:
\begin{itemize}
\item A framework is presented to control a cluster of DHWHs using a \textit{aggregate-and-dispatch} strategy in combination dispatch strategy using reinforcement learning.
\item It is detailed how fitted q-iteration is used to calculate an advantage function for each DHWH which can be used to represent a completely model-free bid-function in a transactive energy setting.
\item The resulting control strategy is applied in hardware in the loop simulation comprising a cluster of actual and simulated domestic hot water heaters.
\end{itemize}

%%%%%%%%%%%%%%%%%%%%%%%%%%%%%%%%%%%%%%
%\vspace{-0.10cm}
%\input{Motivation}
\vspace{-0.10cm}

\section{Implementation}
\label{ImplementationDDR}

\begin{figure}[t!]
\centering{\includegraphics[width=1.0\columnwidth]{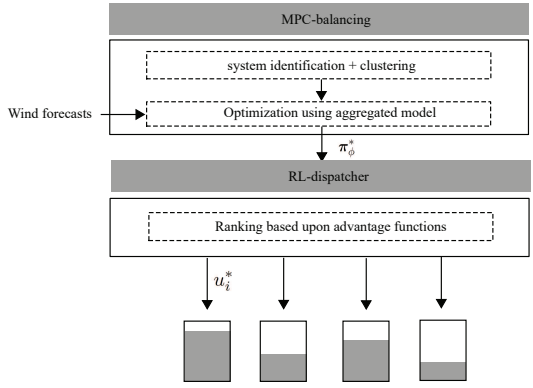}}
\caption{\small Overview of the controller architecture as presented in this paper, data from the domestic hot water systems is used to create a set of aggregated models in an MPC controller that calculates a set point $\pi^{*}_{\phi}$ for a cluster of hot water storage systems. This aggregated set-point is dispatched to the individual assets based upon their individual advantage functions which are learned from data directly in a model-free way.}
\label{fig.controllerOverview}
\end{figure} 

In this Section we present the control framework featuring the fusion of model-based and model-free control strategies. An overview is presented in Fig. \ref{fig.controllerOverview}. This novel control framework comprises two layers: 
\begin{itemize}
\item{In the upper layer, an MPC controller calculates aggregated energy set-points for clusters of DHWHs.}
\item{In the lower layer, a dispatch strategy decomposes the aggregated set-points using advantage functions learned for each individual DHWH.}
\end{itemize}
The framework can be used for providing services in which an entity seeks to control real power consumption from the total aggregation of loads, such as frequency regulation or peak demand shaving; in Section~\ref{sec_simulation} we will explore an application in which the aggregation absorbs wind generator forecast error.

\subsection{MPC control}\label{sec.MPC}
Given the clustered hot water draw profiles, the individual discrete-time dynamics of the water temperature is represented as \cite{Vrettos_ICSGC_2013, Paull_EPSR_2010}
\begin{equation} \label{temperature_dynamics}
\theta_{k+1}^{\imath,\jmath}=a_{k}^{\imath}\theta_{k}^{\imath,\jmath}+\bar{a}_{k}^{\imath}\zeta_{k}^{\imath}+\bar{a}_{k}^{\imath}b_{k}^{\imath}g_{k}^{\imath,\jmath},
\end{equation}
where
\begin{equation}
\begin{aligned}
a_{k}^{\imath}&=1-\bar{a}_{k}^{\imath}=e^{-\frac{\Delta t}{C_{th}}(G+B_{k}^{\imath})}, \\
%\bar{a}_{k}^{\imath}&=1-a_{k}^{\imath}, \\
\zeta_{k}^{\imath}&={G}/({G+B_{k}^{\imath}})T_{k}^a+{B_{k}^{\imath}}({G+B_{k}^{\imath}})T^{in} , \\
%\zeta_{k}^{\imath}&=\frac{G}{G+B_{k}^{\imath}}T_{k}^a+\frac{B_{k}^{\imath}}{G+B_{k}^{\imath}}T^{in} , \\
b_{k}^{\imath}&={\eta}({G+B_{k}^{\imath}}). \nonumber
%b_{k}^{\imath}&=\frac{\eta}{G+B_{k}^{\imath}}.
\end{aligned}
\end{equation}
Herein, $k$ denotes the 1-min discrete time index; $\imath\in\{1,\ldots,n\}$ denotes the cluster number; $\jmath\in\{1,\ldots,n_\imath\}$ denotes the $\jmath$th boiler in the $\imath$th cluster; $\theta_{k}^{\imath,\jmath}$ is the water temperature; $\Delta t$ is the discrete time step; $C_{th}=mc_{th}$, where $m$ is the water mass and $c_{th}$ is the water specific heat; $G=A_l/R_{th}$, where $A_l$ is the total boiler surface and $R_{th}$ is the tank thermal resistance; $B_{k}^{\imath}=\rho_w D_{k}^{w,\imath}c_{th}$, where $\rho_w$ is the water density and $D_{k}^{w,\imath}$ is the $\imath$th cluster's mean hot water usage; $T_{k}^a$ and $T^{in}$ are the ambient temperature and the inlet water temperature, respectively; $\eta$ is the heating efficiency; $g_{k}^{\imath,\jmath}$ is the heating power.

Let $K$ denote the total 1-min time steps in a 15-min window. Augmenting $g_{k}^{\imath,\jmath}$ and $\zeta_{k}^{\imath}$ yields
\begin{equation}
\begin{aligned}
\hat{\boldsymbol{\zeta}}_k^{\imath}&=\left[ \zeta_{k}^{\imath}~\zeta_{k+1}^{\imath}~\ldots ~\zeta_{k+K-1}^{\imath} \right]^{\mathsf{T}}, \\
\hat{\boldsymbol{g}}_{k}^{\imath,\jmath}&=\left[ g_{k|k}^{\imath,\jmath}~g_{k+1|k}^{\imath,\jmath}~\ldots~ g_{k+K-1|k}^{\imath,\jmath} \right]^{\mathsf{T}}. \nonumber
\end{aligned}
\end{equation}
Therefore, the final state can be calculated as
\begin{equation}\label{10ssystem}
\theta_{k+K|k}^{\imath,\jmath}=M_{k}^{\imath,K}\theta_{k}^{\imath,\jmath}+\boldsymbol{C}_{k}^{\imath,K}\hat{\boldsymbol{\zeta}}_k^{\imath}+\boldsymbol{D}_{k}^{\imath,K}\hat{\boldsymbol{g}}_{k}^{\imath,\jmath},
\end{equation}
where
\begin{equation}
\begin{aligned}
M_{k}^{\imath,K}&=\prod_{\kappa=k}^{k+K-1}a_i(\kappa), \\
\boldsymbol{C}_{k}^{\imath,K}&= \left[ \prod_{\kappa=k+1}^{k+K-1}a_i(\kappa)\bar{a}_{k}^{\imath} \prod_{\kappa=k+2}^{k+K-1}a_i(\kappa)\bar{a}_{k+1}^{\imath}  ~ \ldots~ \bar{a}_{k+K-1}^{\imath} \right] , \\
\boldsymbol{D}_{k}^{\imath,K}&=\boldsymbol{C}_{k}^{\imath,K} \circ \left[ b_{k}^{\imath} ~ b_{k+1}^{\imath}~\ldots~b_{k+K-1}^{\imath} \right]. \nonumber
\end{aligned}
\end{equation}

Let $\phi$ denote the 15-min time index. Summing dynamics in \eqref{10ssystem} over $\jmath$ yields the aggregate dynamics
\begin{equation}\label{aggregated_temperature}
\Theta_{\phi+1|\phi}^{\imath}=M_{\phi}^{\imath,K}\Theta_{\phi}^{\imath}+n_\imath \boldsymbol{C}_{\phi}^{\imath,K}\hat{\boldsymbol{\zeta}}_\phi^\imath+\frac{1}{\Delta T}\boldsymbol{D}_{\phi}^{\imath,K} \bold{1} \pi_{\phi}^{\imath},
\end{equation}
where $\Delta T$ is the 15-min time interval, $\Theta_{\phi}^{\imath}=\sum_{\jmath=1}^{n_\imath }\theta_{\phi}^{\imath}$ is the aggregated water temperature, and $\pi_{\phi}^{\imath}$ [J] is the energy setpoint. 
Augmenting system \eqref{aggregated_temperature} by collecting all clusters, we define
\begin{equation}
\begin{aligned}
\boldsymbol{\Theta}_\phi&=\left[ \Theta_{\phi}^1~\Theta_{\phi}^2~\ldots~\Theta_{\phi}^n \right]^{\mathsf{T}} \in \mathds{R}^{n}, \\
\hat{\boldsymbol{\zeta}}_{\phi} &= \left[ \hat{\zeta}_{\phi}^1 ~ \hat{\zeta}_{\phi}^2~ \ldots~\hat{\zeta}_{\phi}^n \right]^{\mathsf{T}} \in \mathds{R}^{nK}, \\
\boldsymbol{\pi}_{\phi}&= \left[ \pi_{\phi}^1~\pi_{\phi}^2~\ldots~\pi_{\phi}^n\right]^{\mathsf{T}} \in \mathds{R}^{n}.
\end{aligned}
\end{equation}
Then in the MPC form we can readily have
\begin{equation} \label{energy_point_system}
\boldsymbol{\Theta}_{\phi+1|\phi}=\boldsymbol{M}_{\phi}^K\boldsymbol{\Theta}_\phi+\boldsymbol{C}_{\phi}^K\hat{\boldsymbol{\zeta}}_{\phi}+\boldsymbol{D}_{\phi}^K\boldsymbol{\pi}_{\phi},
\end{equation}
where $\boldsymbol{M}_{\phi}^K\in \mathds{R}^{n\times n}$, $\boldsymbol{C}_{\phi}^K\in \mathds{R}^{n\times nK}$, and $\boldsymbol{D}_{\phi}^K\in \mathds{R}^{n\times n}$ are block diagonal matrices comprised of $M_{\phi}^{\imath,K}$, $n_\imath \boldsymbol{C}_{\phi}^{\imath,K}$, and $\boldsymbol{D}_{\phi}^{\imath,K}\bold{1}/\Delta T$, respectively, $\imath=1,\ldots,n$. The constant power during $[\phi,\phi+1]$ can be represented as
\begin{equation}
P_{\phi}^b={\bold{1}_n^{\mathsf{T}}\boldsymbol{\pi}_{\phi}}/{\Delta T}.
\end{equation}

%In this work, we aim at making use of the aggregated power consumption to balance the supply-demand via compensating for the deviation of short-term forecasted wind power generation from the long-term planned wind power generation. Specifically,

Backup controllers introduce binary conditions. Let $[\underline{T},\overline{T}]$ denote the targeted temperature range. At time $\phi$, define
\begin{equation}
\boldsymbol{\sigma}_{\phi} =\left[ \sigma_{\phi}^1~\sigma_{\phi}^2~\cdots~\sigma_{\phi}^n \right]^{\mathsf{T}} \in \mathds{B}^{n}, 
\end{equation}
where each scalar entry $\sigma_{\phi}^{\imath}$ is a binary variable defining
\begin{equation} \label{binary_decision}
\left\{ \begin{array}{cc}
\Theta_{\phi}^{\imath}>n_\imath \underline{T}~&\text{if}~\sigma_{\phi}^{\imath}=1, \\
\Theta_{\phi}^{\imath}<n_\imath \underline{T}~&\text{if}~\sigma_{\phi}^{\imath}=0. 
\end{array}\right.
\end{equation}
Since $\left| \Theta_{\phi}^{\imath}-n_\imath \underline{T} \right|$ is always physically bounded, by using the big-$M$ method, we have \eqref{binary_decision} rewritten as
\begin{equation} \label{binary_decision_M}
\begin{aligned}
\Theta_{\phi+\varphi|\phi}^{\imath}&\geq n_\imath\underline{T}+\epsilon-\mathcal{M}_\Theta(1-\sigma_{\phi+\varphi|\phi}^{\imath}), \\
n_\imath \underline{T}~~~~&\geq \Theta_{\phi+\varphi|\phi}^{\imath} +\epsilon - \mathcal{M}_\Theta\sigma_{\phi+\varphi|\phi}^{\imath},
\end{aligned}
\end{equation}
where $\mathcal{M}_{\Theta}\gg \sup\left\{ \left| \Theta_{\phi+\varphi|\phi}^{\imath}-n_\imath \underline{T} \right| \right\}$ and $\epsilon$ is an infinitesimal.  %$\forall i=1,\ldots,n, \iota=0,\ldots, L-1$. 

If the cluster average temperature drops below $\underline{T}$, all boilers in it or at least a certain fraction must be turned on, yielding
\begin{equation} \label{logic_constraint_on_H}
\begin{array}{cc}
{0} \leq \pi_{\phi+\varphi|\phi}^{\imath} \leq \Delta T \bar{P} n_\imath , &~\text{if}~\sigma_{\phi+\varphi}^{\imath}=1, \\
\pi_{\phi+\varphi|\phi}^{\imath} = \Delta T \bar{P} n_\imath , &~\text{if}~\sigma_{\phi+\varphi}^{\imath}=0,
\end{array}
\end{equation}
$\forall~\varphi=0,\cdots,\Phi-1$ and $\Phi$ is the prediction horizon. Note that $\sigma_{\phi|\phi}^{\imath}$ is known, while $\sigma_{\phi+1|\phi}^{\imath},\ldots,\sigma_{\phi+\Phi-1|\phi}^{\imath}$ are decision variables. Equivalently, \eqref{logic_constraint_on_H} can be rewritten as
\begin{equation} \label{inequality_constraint_on_H}
\begin{aligned}
&(1-\sigma_{\phi+\varphi|\phi}^{\imath}) \Delta T \bar{P} n_\imath  \leq \pi_{\phi+\varphi|\phi}^{\imath} \leq \Delta T \bar{P} n_\imath , \\
&\forall~\imath=1,\ldots,n,~\varphi=0,\ldots,\Phi-1.
\end{aligned}
\end{equation}

To ensure an acceptable temperature range, a temperature hard bound $[\underline{\underline{T}},\overline{T}]$, where $\underline{\underline{T}}<\underline{T}$, must be imposed as 
\begin{equation}\label{hard_temp_constraint}
\underline{\underline{T}}n_\imath  \leq \Theta_{\phi+\varphi|\phi}^{\imath} \leq \overline{T}n_\imath , ~\forall \imath=1,\ldots,n, \varphi=1,\ldots,\Phi.
\end{equation}

Therefore, one of the control constraint sets $\mathbb{\Pi}$ is defined as
\begin{equation}
\mathbb{\Pi}:=\left\{ \boldsymbol{\pi}_{\phi+\varphi|\phi} | ~\eqref{binary_decision_M}, \eqref{inequality_constraint_on_H},~\text{and~}\eqref{hard_temp_constraint}~\text{hold} \right\}.
\end{equation}

At time $\phi$, let $\bar{P}_{\phi}^{w}$ and $P_{\phi}^{w}$ denote the planed and short-term forecasted wind power generation, respectively, and $\bar{P}_{\phi}^b$ denote the baseline aggregated boiler power consumption. System balancing within the prediction horizon $\Phi$ can be achieved by solving an MIQP problem as
\begin{equation}
\begin{aligned}
\boldsymbol{\pi}_\phi^* &=\arg\min_{\scriptsize \begin{aligned} &\boldsymbol{\pi}_{\varphi_\pi|\phi} \\ & \boldsymbol{\sigma}_{\varphi_\sigma|\phi} \end{aligned}} \sum_{\varphi=\phi}^{\phi+\Phi-1}\left( P_{\varphi|\phi}^b-\bar{P}_{\varphi}^b-  P_{\varphi}^w+\bar{P}_{\varphi}^w \right)^2  \\
\text{s.t.~}& \boldsymbol{\pi}_{\varphi_\pi|\phi} \in \mathbb{\Pi},~\forall~\varphi_\pi=\phi,\ldots,\phi+\Phi-1, \\
& \boldsymbol{\sigma}_{\varphi_\sigma|\phi} \in \mathds{B}^{n\Phi},~\forall~\varphi_\sigma=\phi+1,\ldots,\phi+\Phi-1.
\end{aligned}
\end{equation}

\subsection{Dispatch}
\label{subsec:dispatch}
\subsubsection{Markov Decision Problem}
As in \cite{convolutional}, the decision making problem is presented as a Markov Decision Process (MDP) defined by the state space $X$, the actions space $U$, the discrete-time transition function $f$ and the cost function $c$. As a result of an action $\mathbf{u}_k$ $\in U$ a state transition occurs from $\x_{k}$ to $\x_{k+1}$ following:
\begin{equation}
\x_{k+1}=f(\x_{k},\action_{k},\w_{k}).
\end{equation}
The random process $\w_{k} \in W$ is drawn from a probability distribution $p_{w}(\cdot,\x_{k})$. The cost function\footnote{In this work, a cost is minimized rather than a reward maximized.} $c$ in turn is defined as:  
\begin{equation}
c_{k}(\x_{k},\action_{k},\x_{k+1})=c_{k}(\x_{k},\action_{k},\w_{k}).
\end{equation} 
The cost function can be defined in a variety of ways.  This flexibility is a central benefit to using reinforcement learning over other demand response dispatch strategies. We will discuss this further when we introduce the concept of an \textit{advantage function} below.

In reinforcement learning a typical objective is to find an optimal state-action value function or $Q$-function that follows the Bellman optimality equation~\cite{BellmanDP}:
\begin{align}
Q^{*}(\x,\action) = \underset{\w\sim\probw}{\mathds{E}}\left[c(\x,\mathbf{u},\w) + \underset{u' \in U}{\text{min~}} Q^{*}(f(\x,\action,\w),\action') \right].
\label{Qfunction}
\end{align}
From this $Q$-function, an optimal policy ${h^{*}:~X~\rightarrow~U}$ is determined by:
\begin{equation}
h^{*}(\x)  \in \underset{\action \in U}{\text{arg min~}} Q^{*}(\x,\action).
\label{Qpolicy}
\end{equation}
This policy can be used to control the domestic hot water heaters.

% Things are a little unclear here. I think this is saying that $h^*$ would be used for local control, without considerations for the aggregation?  Please clarify.  Also need to define $Q^h$.   

As is explained in Section~\ref{subsec.realtime}, our demand response dispatch strategy is to activate (turn on or off) those devices for which the opportunity cost for deviating from a policy $h$ is lowest. We measure that opportunity cost with the advantage function $A^{h}$~\cite{duelling}, defined as: 
\begin{equation}
 A^{h}(\x,\action)=Q^{h}(\x,\action)-V^{h}(\x). 
\label{Qpolicy}
\end{equation} 
Here $V^{h}(\x)$ is the value function:
\begin{align}
V^{h}(\x) &= \underset{\w\sim\probw}{\mathds{E}}\left[c(\x,h(\x),\w) + Q^{h}(\x',h(\x')) \right]\\\nonumber
\x' &=f(\x,\action,\w). 
\label{Qfunction}
\end{align}

\subsubsection{DHWL level implementation}
In this paper the MPD developed by Ruelens \textit{et al}~\cite{RuelensBRLDevice} is used, which we summarize as follows.  For each DHWH with index $i$ in the set $\mathcal{D}$ the state space $X$ comprises: time-dependent state information $\Xtime$ and controllable state information $\Xphys$. The time-dependent component is essential to capture time-dependent patterns such as hot water usage behavior. 
The time-dependent component contains the quarter-hour \footnote{A natural extension is to add the day of the week.} of the day:
\begin{equation}
\xtime \in \Xtime= \left\{1,\dots,96\right\}.
\label{eq.timestate}
\end{equation}

The controllable state information $\xphysk$ is the operational temperature $T_{k}^{i}$ of each DHWH: 
\begin{equation}
\underline{T}_{k}^{i}< T_{k}^{i} <\overline{T}_{k}^{i} \, 
\label{eq.physstate}
\end{equation} 
where $\underline{T}_{k}^{i}$ and $\overline{T}_{k}^{i}$ denote the lower and upper bound set by the end user. %As it will be presented in Sec.~\ref{Sec:controlaction}, $X_{phys}$ will further extended with information on the previous control action.
%%%%%%%%%%%%%%
This results in the observable state vector $\xobski$ for DHWH $i$:
\begin{equation}
\xobski =\left(x_{\mathrm{t}, k},T_{k}^{i}\right). 
\label{eq.stateDef}
\end{equation}

The control action for each DHWH is a binary value indicating if the DHWH is switched ON of OFF:
\begin{equation}
u_{k}^{i} \in \left\{0,1\right\}. 
\label{eq.contr_space}
\end{equation}

As in~\cite{koch2011modeling} and~\cite{RuelensBRLDevice}, we assume each DHWH is equipped with a backup controller, overruling the control action resulting from the policy $h^{i}$. 
Although this can be incorporated directly in the policy, the rationale for separating this is that in a commercial implementation the exact details of the backup controller can be shielded by the manufacturer. 

The function  $B$ maps the requested control action $u_{k}^{i}$ to a physical control action $\uphysk$, depending on its  temperature $T_{k}^{i}$ 
and is defined as:
\begin{align}
\uphysk=B(T_{k}^{i},u_{k}^{i},\chi^{i}) = \left\{\begin{matrix}
1&\text{if }& {T_{k}^{i} \leq}\underline{T}_{k}^{i} \quad \quad \quad \quad \;\;\, \\ 
u_k^{i}&\text{if }& {\underline{T}_{k}^{i} \leq } {T_{k}^{i} \leq}\overline{T}_{k}^{i}.\quad \quad \\ 
0 &\text{if }& {T_{k}^{i} >}\overline{T}_{k}^{i} \quad \quad \quad \quad \;\;\,
\end{matrix}\right.
\label{Eq:backup_controller_2}
\end{align}
Here $\boldsymbol{\chi}^{i}$ contains the minimum and maximum temperature boundaries, $\underline{T}_{k}^{i}$  and $\overline{T}_{k}^{i}$.
The output power of DHWH $i$ when switched on is referred to as $p^{nom,i}$.

The cost function used in this paper comprises two components, one represents the cost for energy and one a availability fee for being available for a demand response event, i.e. to be switched ON, $\alpha$ is the accompanying fee. The cost function $c$  is defined as:
\begin{equation}
c \left(\xobski, \uphysk,\lambda ,\alpha,\right) = p^{nom,i} \Delta t \lambda \uphysk-\alpha I(T_{k}^{i}>\underline{T}_{k}^{i}),
\label{eq.rewardDispatch}
\end{equation}
where $I(T_{k}^{i}>\underline{T}_{k}^{i})$ is the indicator function, equal to 1 where $T_{k}^{i}>\underline{T}_{k}^{i}$ and 0 otherwise, $\Delta t$ corresponds to the length of the activation and $\lambda$ the energy price. 
Note that this cost function can be readily adapted to integrate other objectives such as self-consumption. For example if renewable generation is located under the same grid connection, an extra cost is obtained when exporting this renewable energy, favoring self-consumption.  

\subsubsection{Advantage Function}
\label{subsec.advantage}
In order to decide which DHWHs to switch on we need an advantage function $A^{h,i}$ for each asset in $\mathcal{D}$, this is calculated from $\widehat{Q^{*}}(\x,u)$ calculated from a batch of four tuples $\mathcal{F}$ for each individual DHWH as discussed in \cite{RuelensBRLDevice}. In general a batch $\mathcal{F}$ of four tuples has the following form:
\begin{equation}
\mathcal{F} =\left\{(\x_{l},u_{l},\x_{l}',c_{l}),\, l = 1,...,|\mathcal{F}|\right\},
\label{eq.tuples}
\end{equation}

Algorithm~\ref{forecastedFQI} is an implementation of fitted q-iteration \cite{ernst2005tree} as detailed in \cite{RuelensBRLDevice}. The tuples contain $\x_{l}'$, the successor state to $\x_l$.
An addition to \cite{RuelensBRLDevice} is that a fitted double q-iteration \cite{DoubleFQI} implementation has been used. This is done as $\widehat{Q^{*}}(\x,u)$ is not used only to derive a policy following Eq. (\ref{Qpolicy}), which is impervious to a static bias in the approximation of ${Q^{*}}(\x,u)$. In this work, $\widehat{Q^{*}}(\x,u)$ is also used to calculate $A^{\widehat{h^{*}}}$ to determine the DHWHs that are to be activated taking into account their opportunity cost as expressed in Eq. (\ref{eq.rewardDispatch}). As such the performance of the approach presented is more susceptible to approximation errors $\epsilon = |\widehat{Q^{*}}(\x,u)-{Q^{*}}(\x,u)|$. By working with in a fitted double q iteration setting, this error is reduced.   
After obtaining $\widehat{Q^{*,i}}(\x,u)$ for DHWH $i$, $A^{i}(\x,u)$ is defined as:
\begin{equation}
 A^{i}(\x,u)=\widehat{Q^{*,i}}(\x,u)-\underset{u \in U}{\text{arg min~}} \widehat{Q^{*,i}}(\x,u). 
\label{Qpolicy}
\end{equation}
To approximate the Q- function, an ensemble of extremely randomized trees~\cite{ernst2005tree} was used. Future research is directed towards using more advanced regression architectures, specifically \cite{duelling} a dueling architecture as it allows for a direct representation of the advantage function $A$.   
\begin{algorithm}[t]
\caption{Fitted Q-iteration as detailed in \cite{RuelensBRLDevice}.}
\label{forecastedFQI}
\begin{algorithmic}[1] 
\algsetup{linenosize=\tiny}
\renewcommand{\algorithmicrequire}{\textbf{Input:}}
\REQUIRE $\mathcal{F}=\{\x_{l}, u_l, \x_{l}', \uphysl\}_{l=1}^{\#\mathcal{F}}, \bm{\lambda},\bm{\alpha}$ \\
\STATE let $\widehat{Q}_{0}$ be zero everywhere on $X$ $\times$ $U$
\FOR {$N=1,\ldots,T$}
\FOR {$l = 1,\ldots,\#\mathcal{F}$}
\STATE $~~c_{l} \leftarrow \rho (\x_{l}, \uphysl, \lambda_{l},\alpha_{l})$
\STATE $~~Q_{N,l}\leftarrow c_{l} +\underset{u \in U}{\text{min~}}Q_{N-1}(\x_{l}',u)  $  
\ENDFOR
\STATE use extra trees \cite{geurts2006extremely} to obtain $\widehat{Q}_{N}$ from $\mathcal{T} = \left\{\left((\x_{l},u_{l}),Q_{N,l}\right),l =1,\ldots,\#\mathcal{F}\right\}$
\ENDFOR
\ENSURE $\widehat{Q^{*}}=\widehat{Q}_{N}$
\end{algorithmic}
\end{algorithm}
\subsubsection{Real time control}
\label{subsec.realtime}
Finally as illustrated in Fig. \ref{fig.controllerOverview}, the control action $\pi^{*}_{\phi}$ resulting from the MPC controller is to be dispatched over the cluster of DHWHs. 
This comes down to solving the following optimization problem:

\begin{align}
u_{1}^{*},\ldots,u_{|\mathcal{D}|}^{*}=&\underset{u_{1},\ldots,u_{|\mathcal{D}|}}{\text{arg min~}} \sum_{i=1}^{|\mathcal{D}|}{A^{i}(\x,u)}\\
\textrm{s.t.:~}\sum_{i=1}^{|\mathcal{D}|}u^{i}>&\pi^{*}_{\phi}/\Delta T.
\label{Eq:Dispatch}
\end{align}
An approximate solution\footnote{From the assumption that $A>0$ and an equal nominal power for each DHWH.} to which is found by following Algorithm \ref{dispatchAlgorithm}. 
The relation with transactive energy is that the advantage functions are considered to be bid-functions and (\ref{Eq:Dispatch}) is regarded as a market clearing problem.
\begin{algorithm}[t]
\caption{Dispatch the control set point $\pi^{*}_{\phi}$.}
\label{dispatchAlgorithm}
\begin{algorithmic}[1] 
\algsetup{linenosize=\tiny}
\renewcommand{\algorithmicrequire}{\textbf{Input:}}
\REQUIRE $ \pi^{*}_{\phi},\x_{i},\Delta T, A_{i}(\x,u)~ \forall i \in \mathcal{D}$ \\
\STATE $\mathcal{D}' = \{i \in \mathcal{D} ~|~ T_{k}^{i}<\underline{T}_{k}^{i}\}$  
\STATE $p_{\textrm{min}} = \sum_{i\in \mathcal{D}'}{p_{i}^{\mathrm{nom}}}$
\STATE $p^{D} \leftarrow p_{\textrm{min}}$  
\STATE $\mathcal{V}=\emptyset$  
\FOR {$N= 1,\ldots, |\mathcal{D}|$}
\STATE $i^{*} \in \underset{i \in \mathcal{D} \setminus \mathcal{V}}{\text{arg min~}} A^{i}(\x,1) $ 
\IF{$p^{D}<\pi^{*}_{\phi}/\Delta T$}
\STATE $u^{i^{*}}=1$
\STATE $p^{D} \leftarrow p^{D}+p_{i}^{\mathrm{nom}}$ 
\ELSE
\STATE $u^{i^{*}}=0$
\ENDIF
\STATE $\mathcal{V}\leftarrow \mathcal{V}\cup i^{*}$
\ENDFOR
\ENSURE $u^{1},\ldots,u^{|\mathcal{D}|}$
\end{algorithmic}
\end{algorithm}

%\input{MPC_formulation_20171113}
%\vspace{-0.10cm}
%\input{ProblemFormulation}
%\vspace{-0.10cm}
\vspace{-0.10cm}
\section{Simulation Results \& Analyses} \label{sec_simulation}

The control framework developed in Section \ref{ImplementationDDR} is verified through a networked simulation platform. The MPC and a database for information exchange are set up in Berkeley, USA; the RL dispatcher and simulators are located in Antwerp, Belgium. The scenario we test is the one in which DHWH are dispatched to absorb wind plant forecast errors.  We assume that a forecast for wind plants, $\bar{P}_w$, is produced one day ahead for use in day ahead energy markets, and that a second wind forecast, $P_w$, is produced 15 minutes ahead for use in real time energy markets.  The specific use case we test is the one in which a wind plant seeks to minimize its exposure to real-time energy market price fluctuations by dispatching the DHWH aggregation to track ${P}_w-\bar{P}_w$, the difference between real time and day-ahead wind forecasts. We kept the simulation running for a 40-day period, from which two consecutive days were randomly picked to present the results. 

\begin{figure}[!htb] \centering
\includegraphics[width=0.43\textwidth, trim = 30mm 1mm 29mm 5mm, clip]{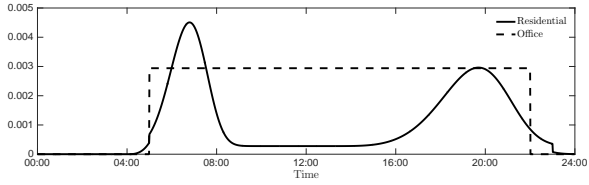}
\caption{Probability distributions of water use events.}
\label{distribution_res_off}
\end{figure}

\begin{figure}[!htb] \centering
\includegraphics[width=0.43\textwidth, trim = 25mm 0mm 25mm 0mm, clip]{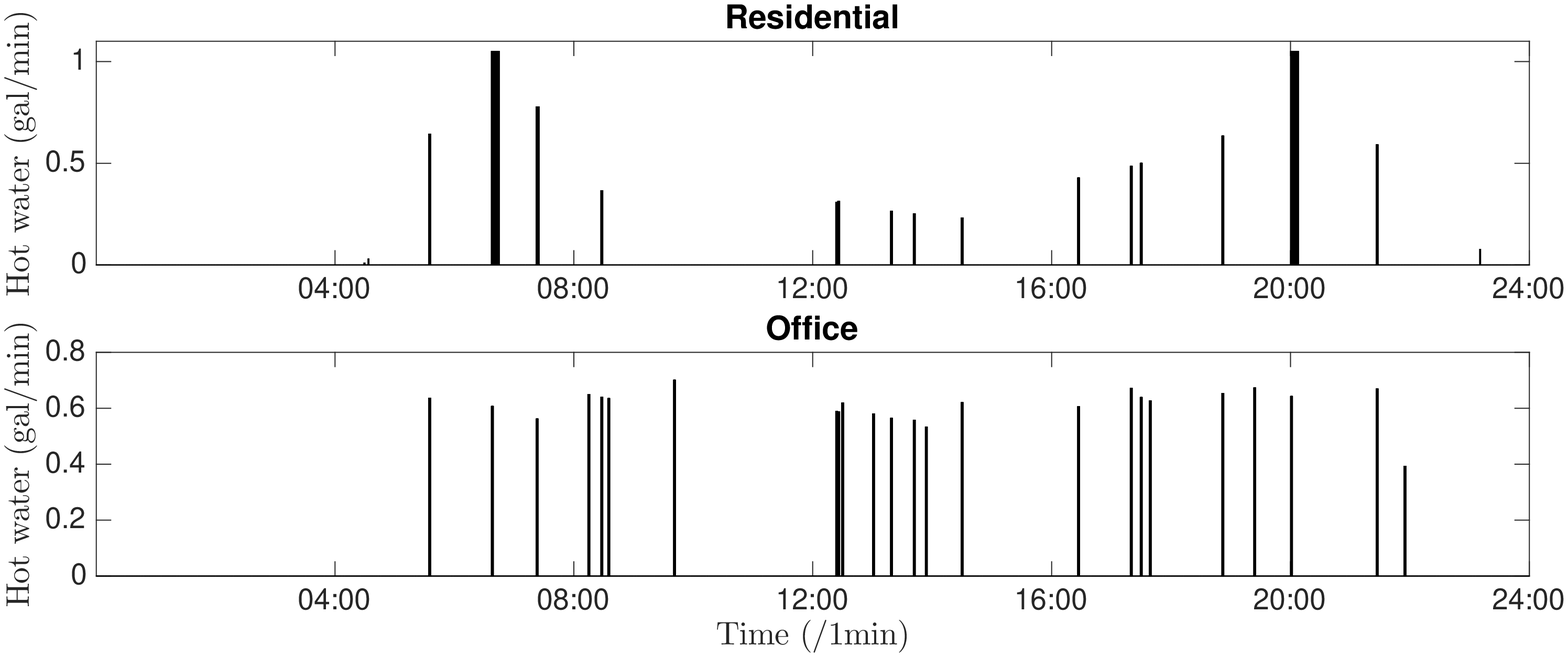}
\caption{Hot water draw profiles.}
\label{hot_water_draw_profiles}
\end{figure}

To set up the simulation, we considered two clusters -- residential and office, and configured 100 simulators for each of them. Probability distributions of hot water use events in the two clusters are shown in Fig. \ref{distribution_res_off}.
We first used the probability distributions to generate 2000 customers' profiles with fixed shower and tap water draw amounts, then calculated the average profile of each cluster, resampled and rescaled it to mimic the real 1-min hot water draw profiles. The mean 1-min hot water draw profiles, which were used for training the RL dispatchers, are shown in Fig. \ref{hot_water_draw_profiles}. It can be readily verified that residential hot water usage has peaks in the morning and night, while office hot water usage is more evenly distributed. The resampled and rescaled profiles were used to generate the baseline DHWH power consumption profile $\bar{P}_b$ via thermostat control. 

We assume that the DHWH aggregation response is measured as the difference between its baseline forecast, $\bar{P}_b$ and its real time consumption produced by Section~\ref{subsec:dispatch}'s RL dispatcher, $P_{b,RL}$.  The RL dispatcher in turn bases its control on dispatch setpoints generated by the MPC controller, ${P}_{b,MPC}$, detailed in Section~\ref{sec.MPC}.

\begin{figure}[!htb] \centering
\includegraphics[width=0.46\textwidth, trim = 20mm 0mm 25mm 0mm, clip]{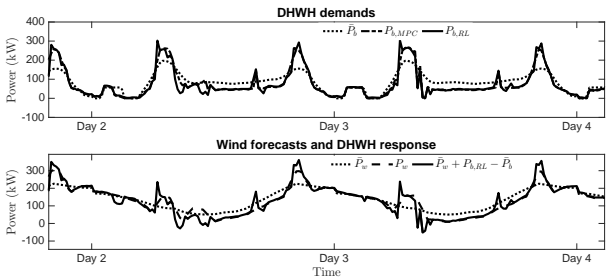}
\caption{All power trajectories.}
\label{wind_power_boiler_power}
\end{figure} 

Fig. \ref{wind_power_boiler_power} shows various power trajectories during the randomly picked two days. The top panel presents three measures of the DHWH demand: baseline, MPC-planned, and RL-dispatched. The bottom panel shows trajectories of day-ahead wind forecast, 15-minute ahead wind forecast, and the sum of day ahead wind forecast and the DHWH response, i.e., $\bar{P}_w + {P}_{b,RL}-\bar{P}_b$.

\begin{figure}[!htb] \centering
\includegraphics[width=0.48\textwidth, trim = 20mm 0mm 25mm 0mm, clip]{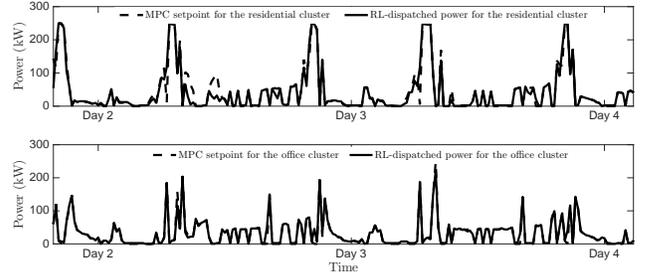}
\caption{Comparisons between MPC and RL for residential and office clusters.}
\label{MPC_RL_comp}
\end{figure}

It can be readily revealed from the top panel of Fig. \ref{wind_power_boiler_power} that both the MPC-planned and RL-dispatched demand trajectories considerably differ from the baseline in order to compensate for wind forecast errors. However one can also see that there are small deviations between the MPC-planned and RL-dispatched trajectories. The percentage normed deviation, $\left\|P_{b,RL}-P_{b,MPC} \right\| / \left\|P_{b,MPC} \right\|$, is 18.29\%. Fig. \ref{MPC_RL_comp} breaks down the deviations into two clusters. The office cluster has a good match, while the residential cluster has few mismatches, leading to relatively big errors in wind power tracking as shown in the bottom panel of Fig. \ref{wind_power_boiler_power}. These mismatches are concentrated in periods of large and long hot water draws, where average temperature inevitably drops below the targeted lower bound.

\begin{figure}[!htb] \centering
\includegraphics[width=0.45\textwidth, trim = 25mm 2mm 33mm 5mm, clip]{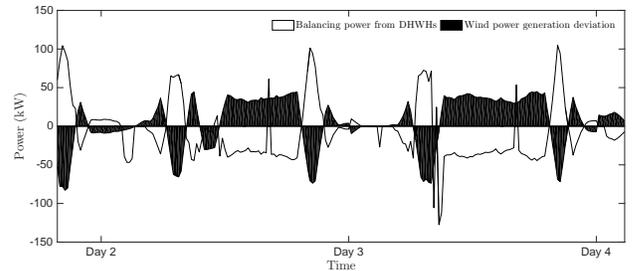}
\caption{Supply-demand balancing performance of MPC-RL.}
\label{supply_demand_balance}
\end{figure}

The bottom panel of Fig. \ref{wind_power_boiler_power} reveals that DHWH response under the proposed control framework can well track the wind power generation, though with a percentage normed deviation of 19.90\%. From another perspective, Fig. \ref{supply_demand_balance} shows the wind power forecast error and the RL-dispatched demand deviation from baseline, both from the supply side. The mean absolute error (MAE) between the RL-dispatched balancing power and wind power generation deviation is about 7.6 kW. 

\begin{figure}[!htb] \centering
\includegraphics[width=0.45\textwidth, trim = 25mm 2mm 33mm 7mm, clip]{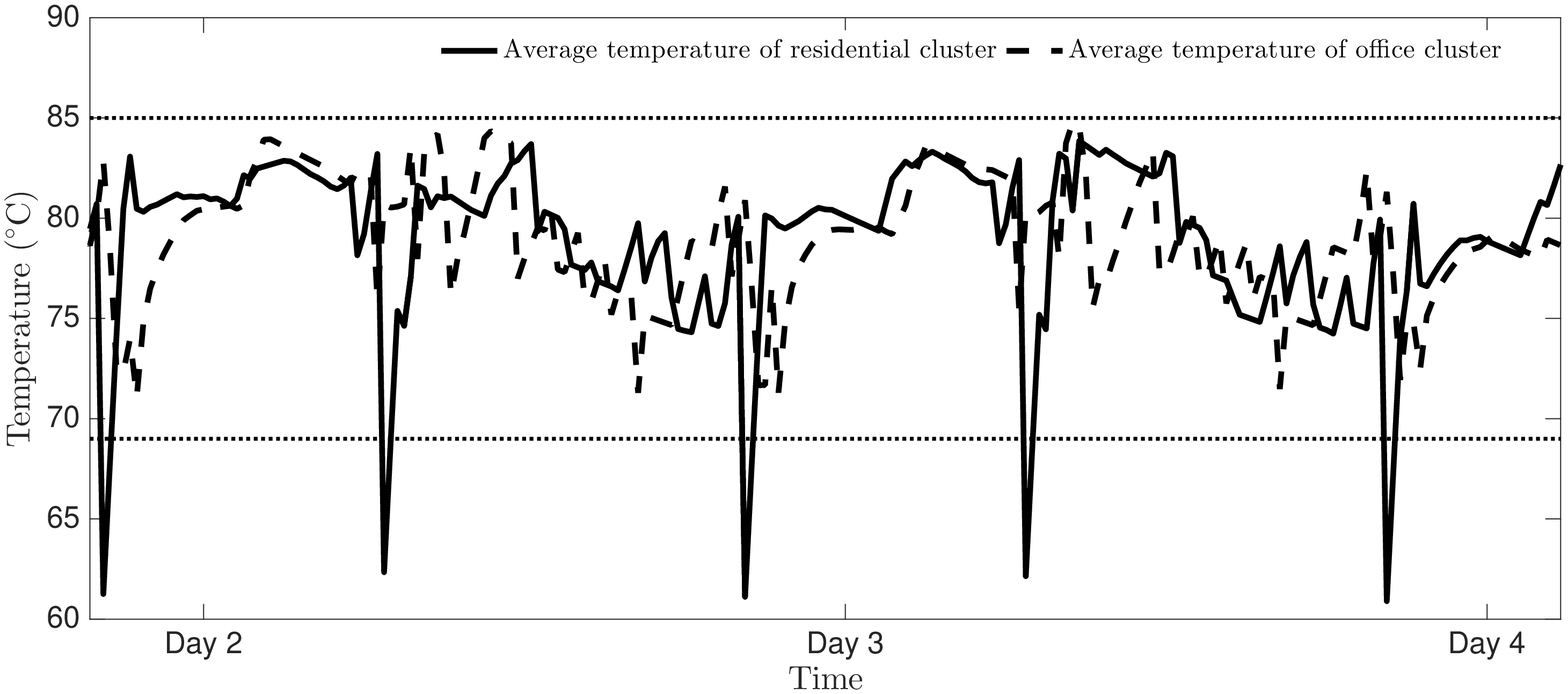}
\caption{Average temperature of residential and office clusters.}
\label{average_temp}
\end{figure}

Under the proposed control framework, average temperature, which was used to calculate MPC setpoints, across DHWHs of two clusters is shown in Fig. \ref{average_temp}.
 Since hard upper bound on temperature is imposed, no temperature exceeds 85$^\circ$C. Due to massive hot water use events, such as shower in the morning and evening, average water temperature of residential DHWHs drops below the targeted lower bound twice a day. In these cases, backup controllers take control and force heating elements of all DHWHs to be in ON status. Despite the fact that the proposed MIQP MPC controller incorporates the backup controller into the setpoint design, the usage of averaged temperature is the main cause for deviations between MPC-planned and RL-dispatched trajectories, as RL dispatchers consider individual temperature which is not necessarily out of the targeted bounds. This phenomenon can be formulated as model uncertainties and will be incorporated into the MPC design in our future work.

\vspace{-0.10cm}
\section{Experiment}
\label{Sec.experiment}
This Sections describes how the approach presented in Section \ref{ImplementationDDR} is deployed in a commercial deployment. 
\subsection{Experimental setup}
The cluster comprised nine General Packet Radio Service (GPRS) connected DHWHs.  Each DHWH contains 80 liters and is equipped with two temperature sensors, a heating element with a power rating of 2.5 kW and a local backup controller that guarantees comfort and safety limits \footnote{The DHWHs are installed at users premises and not in a lab environment. This restricts the duration of the experiment as it can have an effective impact on the energy bill of the end-user as such the experiment has been limited to several hours.}.

Local measurements comprising temperature and power are sent over GPRS to a cloud-based IoT platform, and from there they are pushed to the controller. Here this data is used to determine both the MPC-setpoint $\pi^{*}_{\phi}$ and the individual control actions $u_{i}$ of each DHWH according to Algorithm 2.
We observed the time between sending an activation command and actually having the power available to be in the range of 2-4 seconds.  This suggests the approach could be used for  ancillary services such as frequency response.  Historical data collected from the installed DHWHs over an 18 day period was used to train the RL-dispatcher and the MPC controller model.  

In the experiment, the MPC controller generated a new setpoint every 15 minutes. During each of these intervals, the RL-dispatcher dispatched the DHWHs on a per minute basis using Algorithm 2. This is illustrated in Fig. \ref{fig.controllerOverview}. Please note that the controller is to be improved in a next iteration, e.g. in the current implementation no feedback is foreseen  that can compensate for missing messages and delayed response. An example of a simple feedback strategy is presented in \cite{ClaessensDHNBRL}, more advanced strategies that forecast the correction dynamically are a topic for future research. 

\begin{figure}[t!]
\centering{\includegraphics[width=1.0\columnwidth]{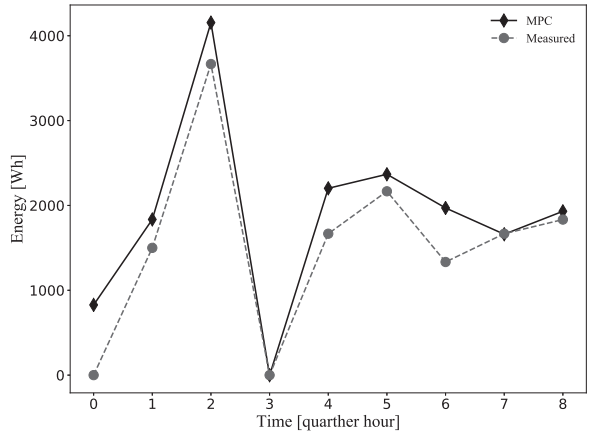}}
\caption{\small Experimental results of the tracking the MPC demand by the cluster of nine RL-dispatched DHWHs. The black diamonds indicate the requested energy over 15 minutes by the MPC controller. The grey circles represent the measured energy over the 15 minutes considered.}
\label{fig.expResult}
\end{figure}

\subsection{Evaluation}
Although the evaluation is limited in scope, it demonstrates the scheme as presented in Section \ref{ImplementationDDR}. For an experiment lasting about 2 hours the MPC controller defined an energy set-point for each quarter hour which was tracked using the RL-dispatcher as presented in Algorithm 2. During the experiment both the advantage function nor the MPC model were updated, only the result from the MPC optimization and the actual values resulting from the advantage functions.  
The results are depicted in Figure \ref{fig.expResult}, the black circles indicate the requested energy from the MPC controller, whilst the grey circles represent the actually observed energy. It is observed that the requested energy is tracked reasonably well, albeit the actual energy consumed is systematically lower. This is due to the fact that some of the commands to switch ON did not reach the DHWH. 

\vspace{-0.10cm}
\section{Conclusions and Future work}
\label{Sec.Conclusions}
This work presents and demonstrates a demand response control architecture for domestic hot water heaters (DHWHs). A model-based control strategy featuring MPC is developed to determine cluster-level control actions using an aggregation model, whilst a model-free dispatch strategy featuring reinforcement learning is used to project the cluster level control action onto DHWH level actions. This results in a scalable and pragmatic control strategy leveraging state of the art in model-based and model-free control which is evaluated using networked simulations. A successful experimental demonstration is provided using a commercial residential demand response implementation.   

In performing this work, several new research questions emerged that will be explored in future work. For example how to use concepts from transfer learning to warm start the advantage function of a DHWH just connected. How to obtain efficient clustering algorithms to cluster the DHWHs, provided the highly uncertain residential water usage profiles. How to design a chance-constrained MPC to accommodate possible temperature violation. How to integrate recent developments in the field of hierarchical learning \cite{hierarchical} to use the DHWS for different applications ranging from ancillary services, self-consumption and local energy services such as peak-shaving.

%\input{Results}
%\vspace{-0.10cm}
%\input{Conclusions}
%\vspace{-0.20cm}

\section{Acknowledgments}
This research is partially supported by H2020 project SIM4BLOCKS, grant number 695965, NSERC, and NSF CyberSEES. The authors would like to thank Itho Daalderop for their support in the experiment. 
%\vspace{-0.20cm}

\ifCLASSOPTIONcaptionsoff
  \newpage
\fi

\bibliographystyle{IEEEtran}  

%\bibliography{C:/Users/ruelensf/Dropbox/bibtext/references} 

\end{document}